\documentclass[11pt]{article}
\usepackage{amsmath,amssymb,amsthm,latexsym}
\usepackage{indentfirst}

\def\0{{\bf 0}}

\theoremstyle{definition}
\newtheorem{theorem}{Theorem}[section]

\newtheorem{lemma}[theorem]{Lemma}

\title{\textbf{An Application of Maximum Principle   to space-like Hypersurfaces with
Constant Mean Curvature in Anti-de Sitter Space}}
\author {{\bf Changxiong Nie } \\
\small{  Faculty of Mathematics and Computer Sciences, Hubei
University, Wuhan, P. R. China}\\
\small{  E-mail: chxnie@163.com}}
\date{}

\begin{document}
\maketitle \footnotetext [ 1]{C. X. Nie is
partially supported by the grant No. 10801006 of NSFC and Zhongdian Natural
Science Foundation of Hubei Educational Committee}
\begin{abstract}
In this paper, we study complete hypersurfaces with constant mean
curvature  in anti-de Sitter space $H^{n+1}_1(-1)$. we prove that if
a complete space-like hypersurface with constant mean curvature
$x:\mathbf M\rightarrow H^{n+1}_1(-1) $ has  two distinct principal
curvatures $\lambda,\mu$, and inf$|\lambda-\mu|>0$, then $x$ is the
standard embedding $  H^{m} (-\frac{1}{r^2})\times   H^{n-m} ( -\frac{1}{1 -  r^2} )$ in anti-de Sitter space $   H^{n+1}_1 (-1)
$.
\end{abstract}



\par\bigskip\noindent
\section{Introduction}
\par\medskip
Let $M^{n+1}_1(c)$ be an (n+1)-dimensional Lorentzian space form of
constant curvature $c$. When $c>0$, $M^{n+1}_1(c)=S^{n+1}_1(c)$,
(n+1)-dimensional de Sitter space; When $c=0$,
$M^{n+1}_1(c)=L^{n+1}$, (n+1)-dimensional Lorentz-Minkowski space;
When $c<0$, $M^{n+1}_1(c)=H^{n+1}_1(c)$, (n+1)-dimensional anti-de
Sitter space. A hypersurface $M$ of $M^{n+1}_1(c)$ is said to be
space-like if the induced metric on $M$ from that of the ambient
space is Riemannian.

The following well-known result of the Bernstein type problem for
maximal space-like hypersurfaces in $M^{n+1}_1(c) (c\geq 0)$ was
proved by Calabi [1], Cheng-Yau [2], and Choquet-Bruhat [6]
\begin{theorem}\cite{ca}\cite{cy}
Let $M$ be an n-dimensional complete maximal space-like hypersurface
in an (n+1)-dimensional Lorentian space form $M^{n+1}_1(c) (c\geq
0)$, then $M$ is totally geodesic.
\end{theorem}
As a generalization of the Bernstein type problem, Cheng-Yau [2] and
T. Ishihara [3] proved that a complete maximal space-like submanifold
$M^n$ of $M^{n+1}_1(c) (c\geq 0)$ is totally geodesic. In [3] T.
Ishihara  also proved the following result
\begin{theorem}\cite{ti}
Let $M$ be an n-dimensional $(n\geq 2)$ complete maximal space-like
hypersurface in anti-de Sitter space $H^{n+1}_1(-1)$, then the norm
square of the second fundamental form of $M$ satisfies
$$S\leq n,$$
and $S=n$ if and only if $M^n=H^m (-\frac{m}{n})\times H^{n-m}
(-\frac{n-m}{n}),(1\leq m\leq n-1)$.
\end{theorem}
In [4], Cao-Wei gave a new characterization of hyperbolic cylinder
$M^n=H^m(-\frac{m}{n})\times H^{n-m} (-\frac{n-m}{n})$ in anti-de
Sitter space $H^{n+1}_1 (-1)$.
\begin{theorem}\cite{cw}
Let $M$ be an n-dimensional $(n\geq 3)$ complete maximal space-like
hypersurface with two distinct principal curvature $\lambda$ and
$\mu$ in anti-de Sitter space $H^{n+1}_1 (-1)$. If
inf$(\lambda-\mu)^2>0$, then $M^n=H^m (-\frac{m}{n})\times H^{n-m}
(-\frac{n-m}{n}),(1\leq m\leq n-1)$.
\end{theorem}
In [4],Cao-Wei also held a conjecture.\\
{\bf Conjecture}: The only complete space-like hypersurfaces in
$M^{n+1}_1(c) (c<0)$ with constant mean curvature and two distinct
principal curvatures $\lambda$ and $\mu$ satisfying
inf$(\lambda-\mu)^2>0$ are the hyperbolic cylinders.

In this paper we investigate complete space-like hypersurfaces in
$M^{n+1}_1(-1)$ with constant mean curvature and two distinct
principal curvatures $\lambda$ and $\mu$ satisfying
inf$(\lambda-\mu)^2>0$, and give an affirmative answer for the
conjecture, and we have the following main theorem.

\begin{theorem}
Let $x:\mathbf M\rightarrow  H^{n+1}_1(-1)$ be an n-dimensional
($n\geq3$) complete   space-like hypersurface in anti-de Sitter
space $H^{n+1}_1(-1) $ with constant mean curvature and with two
distinct principal curvatures $\lambda,\mu$. If
inf$|\lambda-\mu|>0$, then $x$ is the standard embedding $  H^{m} (-\frac{1}{r^2})\times   H^{n-m} ( -\frac{1}{1 -  r^2} )$ in anti-de
Sitter space $   H^{n+1}_1(-1) $.
\end{theorem}
Much recently, however, Wu has more general results like
\begin{theorem}\cite{wu1}
The only complete space-like hypersurfaces in
Lorentz-Minkowski (n+1)-spaces ($n\geq   3$) of nonzero constant $m$th mean curvature
 ($m \leq  n-1$) with two distinct principal
curvatures $\lambda$  and $\mu$ satisfying inf$(\lambda-\mu)^2 > 0$ are the hyperbolic cylinders.
\end{theorem}
We should remind readers that Wu has used Otsuki's idea while we immediately use the maximum principle.
So our proof is more natural and concise. In fact, Wu's results in [10] can be concluded from our method also.
\par\bigskip\noindent
\section{Preliminaries}
\par\medskip
Let $x:\mathbf M  \rightarrow H^{n+1}_1(-1) $ be an n-dimensional
($n\geq3$) space-like hypersurface. Let   $e_1,\cdots, e_n$
 be a local orthonormalbasis of M with respect to the induced metric, and $\omega_1,\cdots,
\omega_n$ their dual form. Let $\xi$ be the local unit normal vector
field such that $\langle \xi,\xi\rangle=-1$.

Denote $x_i=e_i(x)$. Then we have the structure equations

$$dx=  \sum_{i=1}^n\omega_i   x_i,\ dx_i=
\sum_{j=1}^n\omega_{ij}   x_j  +  \mathfrak{h}_i\xi
+\omega_ix,\
 d\xi=  \sum_{i=1}^n \mathfrak{h}_ix_i. \eqno{(2.1)}$$

Denote $\mathfrak{h}_i=\sum_{j=1}^nh_{ij}\omega_j$, from [2] we have
$h_{ij}=h_{ji}$. The curvature tensor can be expressed as Gauss
equation
$$R_{ijkl} = -( h_{ik}h_{jl}-h_{il}h_{jk} )
 -( \delta_{ik}\delta_{jl}-\delta_{il}\delta_{jk} ). \eqno{(2.2)}$$

 And Codazzi equation is
 $$  h_{ijk}= h_{ikj} , \eqno{(2.3)}$$
 where
 $$  \sum_{k=1}^n h_{ijk}\omega_k= dh_{ij}+
 \sum_{k=1}^n h_{kj}\omega_{ki} +
 \sum_{k=1}^n h_{ik}\omega_{kj}. \eqno{(2.4)}$$

The mean curvature of $M$ is given by $H=\frac{1}{n}\sum_ih_{ii}$.
If $H=0$, then $M$ is said to be Maximal, and $H=constant$, then $M$
is said to be of constant mean curvature.

We can choose an appropriate orthonormal basis $e_1,\cdots, e_n$
such that $$h_{ij}=\lambda_i\delta_{ij},$$ where $\lambda_i$ are
principal curvatures.

If we suppose the hypersurface  $x$ has two distinct principal
curvature and has constant mean curvature $H$, then choose an
appropriate orthonormal basis $e_1,\cdots, e_n$ such that
$$ \lambda_1=\cdots=\lambda_m=\lambda,\lambda_{m+1}=\cdots=\lambda_n=\mu,$$
thus we obtain
$$m\lambda+(n-m)\mu=nH=constant, \eqno{(2.5)}$$
{\bf Example 2.1.} Hyperbolic cylinder
$$M_{m,n-m}=
H^{m} (-\frac{1}{r^2})\times   H^{n-m} ( -\frac{1}{1 -  r^2} )
,(1\leq m\leq n-1).$$

We know (see [3]) that $M_{m,n-m}$ is a complete space-like
hypersurface in $H^{n+1}_1(-1)$ with constant mean curvature $H$ and
two distinct principal curvature $\lambda$ and $\mu$, where
$$\lambda_1=\cdots=\lambda_m=\frac{1}{r},\lambda_{m+1}=\cdots=\lambda_n=\frac{1}{\sqrt{1-r^2}}.$$
Thus $M_{m,n-m}$ have constant mean curvature
$H=\frac{m}{nr}+\frac{n-m}{n\sqrt{1-r^2}}$.

Now we have to consider two cases.

Case 1: $2\leq m\leq n-2$.

In this case we make use of the following convention on the ranges
of indices:
  $$1\leq i,j,k \leq m;\ m+1\leq \alpha,\beta,\gamma \leq n
  ;\  1\leq A,B \leq n.$$
{\bf Proposition 2.1.}
 Let $x:\mathbf M\rightarrow  H^{n+1}_1(-1)$ be an n-dimensional
($n\geq3$) complete  space-like hypersurface in anti-de Sitter
space $H^{n+1}_1(-1) $ with constant mean curvature and with two
distinct principal curvatures . If the multiplicities of these two
distinct principal curvatures are greater than one, then $x$ is the
standard embedding $  H^{m} (-\frac{1}{r^2})\times   H^{n-m} ( -\frac{1}{1 -  r^2} )$ in anti-de Sitter space $   H^{n+1}_1
(-1)$.

  {\sl Proof.}
 Letting $i\neq j$ or $\alpha\neq \beta$ in equation(2.4),  there is
 $$ \sum_Ah_{ijA}\omega_A=dh_{ij}+
  \sum_Ah_{ A j}\omega_{Ai} + \sum_Ah_{iA}\omega_{Aj} = \lambda (
  \omega_{ij} + \omega_{ji})= 0,$$
  $$\sum_Ah_{\alpha \beta A}\omega_A=dh_{\alpha \beta}+
  \sum_Ah_{ A \beta}\omega_{A\alpha } + \sum_Ah_{\alpha A}\omega_{A\beta} = \lambda (
  \omega_{\alpha \beta} + \omega_{\beta\alpha })= 0.$$

 That is, when $i\neq j$ and $\alpha\neq \beta$, we have
 $ h_{ijA}= 0, h_{\alpha \beta A} = 0,\forall A.$

By letting $i= j$ or $\alpha= \beta$ in equation(2.4), there is
 $$d\lambda=\sum_{A=1}^n h_{iiA}\omega_A = h_{iii}\omega_i
  + \sum_{\alpha=m+1}^n h_{ii\alpha}\omega_\alpha
 , \forall i,\eqno{(2.6)}$$
 $$ d\mu=\sum_{A=1}^n h_{\alpha\alpha A}\omega_A = h_{\alpha\alpha\alpha}\omega_\alpha
 + \sum_{i=1}^n h_{i\alpha\alpha}\omega_i,  \forall  \alpha. \eqno{(2.7)}$$

 Since $2\leq m \leq n-2$,
equations (2.6) and (2.7) come to
 $$d\lambda= \sum_{\alpha=m+1}^n h_{ii\alpha}\omega_\alpha
 ,   d\mu= \sum_{i=1}^n h_{i\alpha\alpha}\omega_i . \eqno{(2.8)}$$

 Since $m\lambda+(n-m)\mu=nH=constant$, we know that
 $$m d\lambda+(n-m) d\mu = 0 .\eqno{ (2.9)}$$

   Combining with equations (2.8) and (2.9), we get
 $$  \lambda=constant,\ \mu=constant,\   h_{ii\alpha}=0
 , \    h_{i\alpha\alpha}=0, \forall i,\alpha.\eqno{(2.10)}$$
 Then we complete the proof of Lemma 2.1.

For any $i$ and $\alpha$ in equation(2.4), we have
 $$\sum_{A=1}^n h_{i\alpha A}\omega_A = h_{ii\alpha}\omega_i
 + h_{i\alpha\alpha}\omega_\alpha
 =0
 = (\lambda - \mu ) \omega_{i\alpha}. $$

That is,
$$\omega_{i\alpha}=0,\text{and}\ (M,I)=(M_1,I_1)\times(M_2,I_2).$$

 We assert that $M_1,M_2$ have constant curvature.
 For $i\neq j$ and $\alpha\neq \beta$, from equation(2.2)
the sectional curvature of $(M_1,I_1)$ and $(M_2,I_2)$ is
 $$K(e_i\wedge e_j) = R_{ijij} = -1-\lambda^2,
K(e_\alpha\wedge e_\beta) = R_{\alpha  \beta\alpha  \beta} =
  -1-\mu^2, \eqno{(2.11)}$$
  respectively.

On other hands from $K(e_i\wedge e_\alpha) = R_{i\alpha i\alpha} =
-1-\lambda\mu=0.$
 Then we know that when $2\leq m \leq n-2$,
 $x(M)$ is locally Lorentz congruent to the standard
embedding $
H^{m} (-\frac{1}{r^2})\times   H^{n-m} ( -\frac{1}{1 -  r^2} )
\subset  H^{n+1}_1(-1) $.

Thus we complete the proof of proposition.

Case 2: $m=n-1$.

 In this case we make use of the following
convention on the ranges of indices:
  $$1\leq i,j,k \leq n-1; 1\leq A,B \leq n.$$

From (2.5), we can suppose that
 $$(n-1)\lambda+\mu=nH,\ \lambda-\mu=n(\lambda-H)\neq0.\eqno{(2.12)}$$

 Similarly, we have
 $$h_{ijA}=0,\forall A, i\neq j;\ \ d\lambda= h_{iii}\omega_i + h_{iin}\omega_n,
 d\mu= \sum_{i=1}^{n-1}h_{inn}\omega_i + h_{nnn}\omega_n,\forall i.\eqno{(2.13)}$$

 Because $n-1\geq2$,
 from equations  (2.13)  and (2.9), we get
 $$d\lambda=   h_{iin}\omega_n,\forall i;\ \
 d\mu=h_{nnn}\omega_n,\ h_{inn}=0,h_{iii}=0,\forall i,  \eqno{(2.14)}$$

 Equation (2.11) comes to
  $$ \omega_{in} = \frac{1}{\lambda-\mu} h_{iin}\omega_i
  = \frac{1}{n ( \lambda -  H  ) } h_{iin}\omega_i.  \eqno{(2.15 )}$$

 And we assert that the integral curve  of $e_n$ is a geodesic because
 $$\nabla_{e_n}e_n=\sum_{i=1}^{n-1}\omega_{ni}(e_n)e_i
 = - \sum_{i=1}^{n-1}\frac{1}{\lambda-\mu} h_{iin}\omega_i
 (e_n)e_i=0.$$
 We also have $d\omega_n=\sum_{i=1}^{n-1}\omega_{ni}\wedge\omega_i=0
 $.  It means that there exists an arc parametric $s$ of the integral curve
  of $e_n$ such that
 $\omega_n=ds$.  Since $M$ is complete, the arc $s$ tends to infinity.

 If we denote $\dot{f}=\frac{df}{ds}$ for any smooth function
 $f=f(s)$ on the integral curve
  of $e_n$, it follows from equation (2.14) that
 $$d\lambda=   \dot{\lambda} ds,\  h_{iin} = \dot{\lambda},\ \  \forall i.
  \eqno{(2.16 )}$$

From equations (2.15) and (2.16) it follows that
 $$\omega_{in}=\frac{\dot{\lambda}}{n(\lambda-H)}\omega_i,\eqno{(2.17 )}$$

 Exploring into
 $$  d\omega_{in}=d (\frac{\dot{\lambda}}{n(\lambda-H)}\omega_i) =
  \sum_{j=1}^{n-1} \omega_{ij}  \wedge  \omega_{jn} - \frac{1}{2}
   \sum_{A,B=1}^{n }
  R_{inAB}  \omega_{A}  \wedge  \omega_{B} \eqno{( 2.18)}$$
  and collecting the items of $\omega_{i}  \wedge  \omega_{n}$,
  we get
  $$\frac{\ddot{\lambda}}{n(\lambda-H)}-\frac{n+1}{n^2}
  \frac{(\dot{\lambda})^2}{(\lambda-H)^2}=R_{inin}=-1-\lambda\mu.  \eqno{(2.19 )}$$
We introduce the following generalized Liouville-type theorem (see
Choi-Kwon-Sun [5]) in order to prove our main theorem.
\begin{theorem}([5])
Let $M$ be a complete Riemannian manifolds whose Ricci curvature is
bounded from below. Let $F$ be any formula of the variable $f$ with
constant coefficients such that
$$F(f)=c_0f^n+c_1f^{n-1}+\cdots+c_kf^{n-k}+c_{k+1},$$
where $n>1,1\geq n-k\geq 0$ and $c_0>c_{k+1}$. If a
$C^2$-nonnegative function $f$ satisfies
$$\Delta f\geq F(f),$$
then we have
$$F(f_1)\leq 0,$$
where $f_1$ denotes the supermum of the given function.
\end{theorem}

\par\bigskip\noindent
\section{Proof of the main theorem}
\par\medskip
In order to complete the proof of our main theorem, we only consider
Case 2. At first, we prove the following key lemma
\begin{lemma}
 Let $x:\mathbf M\rightarrow  H^{n+1}_1(-1)$ be an n-dimensional
($n\geq3$) complete  space-like hypersurface in anti-de Sitter space $H^{n+1}_1(-1) $ with constant mean curvature and two
distinct principal curvatures . If  one of two principal curvatures
is simple, then Ricci curvature of $M$ is negative semi-definite.
\end{lemma}
{\sl Proof.} From Gauss equation (2.2) and $(n-1)\lambda+\mu=nH$ we
get that
\begin{equation*}
\begin{split}
&Ric_{nn}=(n-1)[\mu^2-nH\mu-1]=(n-1)[(\mu-\frac{n}{2}H)^2-(\frac{n^2}{4}H^2+1)],\\
&Ric_{ii}=-(n-1)-nH\lambda+\lambda^2=(\lambda-\frac{n}{2}H)^2-(\frac{n^2}{4}H^2+n-1),1\leq
i\leq n-1.
\end{split}
\end{equation*}
Thus we have
$$Ric_{nn}\geq-(n-1)(\frac{n^2}{4}H^2+1),Ric_{ii}\geq -(\frac{n^2}{4}H^2+n-1),1\leq
i\leq n-1,$$ so Ricci curvature of $M$ is bounded from below.

From (2.14) we have
\begin{equation*}
 \Delta(\lambda-H)=\sum_A(e_Ae_A-\nabla_{e_A} e_A )(\lambda-H)
 =\ddot{\lambda}-\frac{(n-1)\dot{\lambda}}{n(\lambda-H)},
\end{equation*}
from this above formula and (2.19) and $(n-1)\lambda+\mu=nH$ we
obtain
\begin{equation}
\begin{split}
\Delta(\lambda-H)&=\frac{2}{n}(\dot{\lambda})^2\\
&+n(n-1)(\lambda-H)^3+n(n-2)H(\lambda-H)^2-n(1+H^2)(\lambda-H).
\end{split}
\end{equation}
We define the formula of the variable $x$ with constant coefficients
$$F(x)=n(n-1)x^3+n(n-2)Hx^2-n(1+H^2)x.$$
Then $C_0=n(n-1)>C_3=0$. From (3.1) we have
$$\Delta(\lambda-H)=\frac{2}{n}(\dot{\lambda})^2+F(\lambda-H)\geq F(\lambda-H).$$
If necessaries, take $\tilde{\xi}=-\xi$ as local unit normal vector
field of $M$, we can assume that $\lambda-H>0$. So
$$sup(\lambda-H)>0.\eqno{(3.2)}$$

From generalized Liouville-type theorem [5] we have
$$F(sup(\lambda-H))\leq 0.\eqno{(3.3)}$$

From Gauss equation (2.2) and $(n-1)\lambda+\mu=nH$ we get
$$F(\lambda-H)=n(\lambda-H)R_{inin}. \eqno{(3.4)}$$
From (3.2) and (3.4) we obtain
$$(n-1)(sup(\lambda-H))^2+(n-2)Hsup(\lambda-H)-(1+H^2)\leq 0.\eqno{(3.5)}$$
Let $f(x)=(n-1)x^2+(n-2)Hx-(1+H^2)$, then $lim_{x\rightarrow
\infty}f(x)=+\infty$.

Since $\lambda-H>0$ and (3.5) we obtain
$$f(\lambda-H)\leq 0,i.e.,R_{inin}\leq 0.\eqno{(3.6)}$$
Thus from Gauss equation and (3.6) we get that
\begin{equation*}
\begin{split}
&Ric_{ii}=\sum_{j\neq i}R_{ijij}+R_{inin}=-1-\lambda^2+R_{inin}\leq 0,\\
&Ric_{nn}=\sum_iR_{inin}=(n-1)R_{inin}\leq 0.
\end{split}
\end{equation*}
Thus we complete the proof of Lemma 3.1.

If we denote $w = |\lambda-H|^{ -  \frac{1}{n}  }$,
 it follows that
 $$   \ddot{w} + w R_{inin} = 0.  \eqno{(3.7 )}$$

  While
 $$  R_{inin}= \frac{1}{n-1}R_{nn}= -1 -\lambda\mu,\ (n-1)\lambda + \mu =
 nH,\
 \lambda=H\pm w^{-n},$$
 we have
$$  \lambda   \mu = H^2 \pm (2-n)H w^{-n} + (1-n) w^{-2n} ,  \eqno{(3.8 )}$$
$$   \ddot{w} -w[ 1+H^2 \pm (2-n)H w^{-n} + (1-n) w^{-2n} ] = 0.  \eqno{(3.9 )}$$

The left hand side of equation (3.9) multiplied by $2\dot{w}$ is
precisely the derivative of the left hand side of the following
equation
$$   \dot{w}^2 -w^2[ 1+  ( H \pm   w^{-n} )^2 ] = C =constant.  \eqno{(3.10 )}$$

 Since $ \ddot{w}= -w R_{inin} = - \frac{1}{n-1}wR_{nn}$ is positive,
 we know that $ \dot{w}$ is monotone.
 Because inf $|\lambda-\mu|>0$, sup $\{w(s)|-\infty<s<+\infty\}$
 is a bounded number. Then $\lim_{s\rightarrow+\infty}\dot{w} $
 or $\lim_{s\rightarrow -\infty}\dot{w} $ cannot be infinity.
 We assert that
 $\lim_{s\rightarrow\infty}\dot{w}=0$.

 In fact, if we suppose that
 $\lim_{s\rightarrow+\infty}\dot{w}=a>0$,
 then $\lim_{s\rightarrow+\infty} w = +\infty$.
 Therefore we immediately know that  equation (3.10) cannot hold  when $s$ tends to
 infinity. On the other hand,
 if we suppose that
 $\lim_{s\rightarrow+\infty}\dot{w}=a<0$,
 then $\lim_{s\rightarrow+\infty} w = -\infty$.
 But we  know $w>0$, which is a contraction.
 Therefore, $\lim_{s\rightarrow\infty}\dot{w}=0$.
 Adding the monotonicity of $\dot{w}$, it follows that $\dot{w}\equiv0$.
 That is, $\lambda$ is constant, and so as
 $\mu$. Similar to the discuss in case 1,
  we know that when $  m = n-1$,
 $x(M)$ is locally Lorentz congruent to the standard
embedding $ H^{n-1} (-\frac{1}{r^2})\times   H^{1} ( -\frac{1}{1 -  r^2} )
\subset  H^{n+1} $.

Thus we complete the proof of Theorem 1.4.

{\bf Acknowledgements: }The author would like to express
gratitude to Professor Changping   Wang
   for his warm-hearted inspiration.

\end{document}